\documentclass[12pt]{article}
\usepackage{amscd, amssymb, latexsym, amsmath}

\newtheorem{thm}{Theorem}

\newtheorem{lem}{Lemma}
\newtheorem{dfn}{Definition}

\newtheorem{rem}{Remark}

\newenvironment{proof}
{\noindent {\em Proof.}} {\hfill $\Box$}
\newcommand{\f}{\frac}
\newcommand{\no}{\noindent}

\numberwithin{thm}{section} \numberwithin{cor}{section}
\numberwithin{pro}{section}
\numberwithin{lem}{section} \numberwithin{equation}{section}

\newcommand{\R}{\mathbb R}

\numberwithin{equation}{section}

\newcommand{\pai}{\frac{\partial}{\partial x^i}}

\begin{document}
\title{A Note on the Stability and Uniqueness for Solutions to  the Minimal Surface
System}
\author{Yng-Ing Lee* and Mu-Tao Wang**}
\maketitle
\leftline{*Department of Mathematics and Taida Institute of Mathematical Sciences,}
\leftline{\text{ }National Taiwan University, Taipei, Taiwan}
\leftline{\text{ }National Center for Theoretical Sciences, Taipei Office}
\centerline{email: yilee@math.ntu.edu.tw}

\leftline{**Department of Mathematics,
Columbia University,  New York, NY 10027,USA}
\centerline{email: mtwang@math.columbia.edu}

\begin{abstract}
In this note, we show that the solution to the Dirichlet problem
for the minimal surface system in any codimension is unique in the space of
distance-decreasing maps. This follows as a corollary of the
following stability theorem: if a minimal submanifold $\Sigma$ is
the graph of a (strictly) distance-decreasing map, then $\Sigma$
is (strictly) stable. It is known that a minimal graph of codimension one
is stable without assuming the distance-decreasing condition. 
We give another criterion for the stability in terms of the two-Jacobians of the map
which in particular covers the codimension one case. 
All theorems are proved
in the more general setting for minimal maps between Riemannian
manifolds. The complete statements of the results appear in Theorem~3.1,
Theorem~3.2, and Theorem~4.1.
\end{abstract}
%
\section{Introduction}
\label{intro}
Let $\Omega$ be a bounded domain in $\R^n$. Recall a $C^2$
vector-valued function $f=(f^1,\cdots, f^m):\Omega\rightarrow
\R^m$ is said to be a solution to the minimal surface system
(see Osserman \cite {os} or Lawson-Osserman \cite{lo}) if

\begin{equation}\label{min_eq1}
\sum_{i,j=1}^n\pai(\sqrt{g}g^{ij} \frac{\partial
f^\alpha}{\partial x^j})=0\,\,\text{for each } \alpha=1\cdots m
\end{equation} where
$g_{ij}=\delta_{ij}+\sum_\alpha \frac{\partial f^\alpha}{\partial
x^i}\frac{\partial f^\alpha}{\partial x^j}$,  $g=\det g_{ij}$ and $
g^{ij} $ is the $(i,j)$ entry of the inverse matrix of $(g_{ij})$.
The graph of $f$ is called a non-parametric minimal submanifold.
Equation (\ref{min_eq1}) is indeed the Euler-Lagrange equation of
the volume functional $\int_\Omega \sqrt{g}dx^1\wedge\cdots\wedge dx^n$.

 In the codimension one case, i.e. $m=1$,  a
simple calculation shows
$g^{ij}=\delta_{ij}-\frac{f_if_j}{1+|\nabla
f|^2}$ and the equation is equivalent to the familiar one,

\begin{equation}\label{min_eq2}
div(\frac{\nabla f}{\sqrt{1+|\nabla f|^2}})=0.\end{equation}

It is well-known that the solution to (\ref{min_eq2}) subject to
the Dirichlet boundary condition is unique and stable(see for
example, Lawson-Osserman \cite{lo}).

However in the higher codimension case ( $m>1$), Lawson and
Osserman \cite{lo} discover a remarkable counterexample to the
uniqueness and stability of solutions of (\ref{min_eq1}) when
$n=m=2$. They construct two distinct non-parametric minimal
surfaces with the same boundary. Lawson and Osserman then show an
unstable non-parametric minimal surface with the same boundary
exists as a result of the theorems of Morse-Tompkins \cite{mt} and
Shiffman \cite{sh}. In the same paper, Lawson and Osserman 
show the Dirichlet problem for the minimal surface system may not
be solvable in higher codimension.

In this paper, we first derive a stability criterion for the
minimal surface system in higher codimension. To describe the
results, we define distance-decreasing maps.

\begin{dfn}A map $f:\Omega \subset \R^n \rightarrow \R^m$ is
called  distance-decreasing if the differential $df$ satisfies $
|df(v)|\leq |v|$  at each point of $\Omega$ for any nonzero vector $v\in \R^n$. It is called strictly
distance-decreasing if $  |df(v)|< |v|$  at each point of
$\Omega$ for any nonzero vector $v\in \R^n$.
\end{dfn}

We prove the following stability theorem.

 \vskip 10pt \noindent
{\bf Theorem A (see Theorem 3.1)} {\it  Suppose  a nonparametric
minimal submanifold $\Sigma$ is the graph
 of a  distance-decreasing map $f:\Omega\subset \R^{n}\rightarrow \R^{m}$. Then
 $\Sigma$ is  stable. It is strictly stable if $f$ is  strictly
distance-decreasing. } \vskip 10pt

This theorem generalizes the stability criterion in \cite{lw}. It
turns out the volume element is a convex function on the space of
distance-decreasing linear transformations. The convexity is
further exploited to derive a uniqueness criterion. Namely, we
show the solution to the Dirichlet problem for the minimal surface
system is unique in the space of distance-decreasing maps.

\vskip 10pt \noindent {\bf Theorem B (see Theorem 3.2)} {\it
Suppose that
 $\Sigma_{0}$ and $\Sigma_{1}$
 are nonparametric minimal
  submanifolds which are the graph
 of  $f_{0}:\Omega\subset \R^{n}\rightarrow \R^{m}$ and
 $f_{1}:\Omega\subset \R^{n}\rightarrow \R^{m}$ respectively.
 If both $f_{0}$ and $f_{1}$ are  distance-decreasing and $f_{0}=f_{1}$
 on $\partial \Omega$, then $\Sigma_{0}=\Sigma_{1}.$} \vskip
10pt We remark that solutions to the Dirichlet problem of minimal
surface systems in higher codimensions are
constructed in \cite{wa1} and the solutions are graphs of
distance-decreasing maps. For earlier uniqueness theorems for
minimal surfaces, we refer to Meek's paper \cite{me}.

We prove slightly more general stability and uniqueness theorems
for minimal maps between Riemannian manifolds in this paper. It
turns out the only extra assumption is on the sign of the curvature of the target
manifold. In particular, Theorem 3.1 implies Theorem A
while Theorem 3.2 implies Theorem B.

It is well-known that any minimal graph of codimension one is volume-minimizing by a calibration argument. To
connect to the codimension one case, we develop
another stability criterion for the
minimal surface system in any codimension in section~4. The criterion is in terms
of the two-Jacobians of $f$. To describe the results, we first recall some
notations. Let $L:\R^{n}\rightarrow \R^{m}$ be a linear transformation.
It induces a linear transformation $\wedge^{2}L$, from the wedge product
$\wedge^{2} \R^{n}$ to $\wedge^{2} \R^{m}$ by
$$(\wedge^{2}L)(v\wedge w)=L(v)\wedge L(w).$$ With this we define
$$|\wedge^{2}L|=\sup _{|v\wedge w|=1}|(\wedge^{2}L)(v\wedge w)|.$$
In particular, $|\wedge^{2}L|=0$ if $L$ is of rank one.

 \vskip 10pt \noindent
{\bf Theorem C (see Theorem 4.1)} {\it  Suppose a nonparametric
minimal submanifold $\Sigma$ is the graph
 of a  map $f:\Omega\subset \R^{n}\rightarrow \R^{m}$.  Then $\Sigma $ is stable if
  $|\wedge^{2} df|(x)\leq \f{1}{n-1}$}. \vskip 10pt
A more refined and more general version is proved in Theorem 4.1.
 The rank of the defining function $f$ of a nonparametric minimal submanifold of codimension one is at most one and thus $|\wedge^{2} df|(x)=0$.
   We prove the results for
minimal maps between Riemannian manifolds as stated in Theorem 4.1.

The first author visited Columbia University
while the paper was under revison. She would like to thank the math department for its hospitality.
She is partially supported by an NSC grant in Taiwan.
The second author is grateful to Ben Andrews and Brian White for
inspiring discussions. He is partially supported by an NSF grant
and a Sloan fellowship.

\section{A non-parametric variation formula for graphs}
Suppose that $(M, g)$ and $(N, h)$ are two Riemannian manifolds.
We fix a local coordinate system $\{x^i\}$ on $M$. Let $f$ be a
smooth map from $(M,g)$ to $(N,h)$. The graph of $f$ is an embedded
submanifold of the product manifold $M\times N$, the induced
metric is given by
$$ \sum_{i,j=1}^{n} G_{ij}dx^{i}dx^{j}=\sum_{i,j=1}^{n} (g_{ij}+
 \langle \,df(\f{\partial  }{\partial x^{i}}),df(\f{\partial  }{\partial x^{j}})\,\rangle  )
 dx^{i}dx^{j},$$ and the  volume of the graph is

$$A=\int_M \sqrt{\det G_{ij}}dx^1\wedge\cdots\wedge dx^n=\int_M dv.$$
Assume that there is a family of maps $f_{t}$, $0\leq t\leq
\epsilon$ from $M$ to $N$ with $f_{0}=f$ on $M$ and $f_{t}=f$
outside a compact subset of $M$. When the boundary of $M$ is
nonempty, we require that $f_{t}=f$ on $\partial M$. In the
following, we compute the first and second variations of the
volumes of the graphs. The variation of the volume form is 
\[\begin{split}  \f{d\sqrt{\det G_{ij}(t)}}{dt}&=
 \f{1}{2}\sum_{i,j}G^{ij}(t)\dot{G}_{ij}(t) \sqrt{\det G_{ij}(t)},
  \end{split}\]
where $ G^{ij}(t) $ is the $(i,j)$ entry of the inverse matrix of
$(G_{ij}(t))$.

 Denote the variation field $\f{df_{t}}{dt}$ by $V(t)$.
For simplicity, we omit the dependency of $G_{ij}$ and $V$ on $t$
in the following calculation. Then
\[\begin{split} \dot{G}_{ij}&=\langle \,\nabla
_{V}df_{t}(\f{\partial  }{\partial x^{i}}), df_{t}(\f{\partial
}{\partial x^{j}})\,\rangle  +\langle \,df_{t}(\f{\partial
}{\partial x^{i}}),
\nabla _{V}df_{t}(\f{\partial  }{\partial x^{j}})\,\rangle  \\
&=\langle \,\nabla _{df_{t}(\f{\partial  }{\partial x^{i}})}V,
df_{t}(\f{\partial  }{\partial x^{j}})\,\rangle  +\langle
\,df_{t}(\f{\partial  }{\partial x^{i}}), \nabla _{
df_{t}(\f{\partial  }{\partial x^{j}})}V\,\rangle .
\end{split}\] Here $\nabla$ is the Riemannian connection on $N$,
and
 $V$ and $df_t(\f{\partial  }{\partial
x^{i}})$ are vector fields tangent to $N$.

 Hence the
first variation formula is
\begin{equation}\label{first}
\begin{split}\f{d A_{t}}{dt}=\int_{M}\sum_{i,j}G^{ij}
\langle \,\nabla _{df_{t}(\f{\partial  }{\partial x^{i}})}V,
df_{t}(\f{\partial  }{\partial x^{j}})\,\rangle\, dv_{t}.\end{split}
\end{equation}
Continuing the computation,  we derive
\begin{equation}\label{2nd}
\begin{split}\f{d^{2} A_{t}}{dt^{2}}&=\f{1}{2}\int_{M}(\sum_{i,j}G^{ij}\ddot{G}_{ij}
-\sum
_{i,j,k,l}G^{ik}\dot{G}_{kl}G^{lj}\dot{G}_{ij})\,dv_{t}+\f{1}{4}\int_{M}
(\sum_{i,j} G^{ij}\dot{G}_{ij} )^{2}\,dv_{t}.\\
\end{split}
\end{equation}
Now
\[\begin{split} \ddot{G}_{ij} &=\langle \,\nabla _{V}\nabla _{df_{t}(\f{\partial  }{\partial x^{i}})}V,
df_{t}(\f{\partial  }{\partial x^{j}})\,\rangle  +\langle
\,df_{t}(\f{\partial  }{\partial x^{i}}), \nabla _{V}\nabla _{
df_{t}(\f{\partial  }{\partial x^{j}})}V\,\rangle  \\&\hskip.5cm+
2\langle \,\nabla _{V}df_{t}(\f{\partial  }{\partial x^{i}}),
\nabla _{V}df_{t} (\f{\partial  }{\partial x^{j}})\,\rangle
\\&=\langle \,R(V,df_{t}(\f{\partial  }{\partial x^{i}}))V,df_{t}(\f{\partial  }{\partial x^{j}})\,\rangle  +
\langle \,\nabla _{df_{t}(\f{\partial  }{\partial x^{i}})}\nabla
_{V}V, df_{t}(\f{\partial  }{\partial x^{j}})\,\rangle
\\&\hskip.5cm + \langle \,R(V,df_{t}(\f{\partial  }{\partial
x^{j}}))V,df_{t}(\f{\partial  }{\partial x^{i}})\,\rangle +\langle
\,df_{t}(\f{\partial  }{\partial x^{i}}), \nabla _{
df_{t}(\f{\partial  }{\partial x^{j}})}\nabla _{V}V\,\rangle
\\&\hskip.5cm+2\langle \,\nabla _{df_{t}(\f{\partial  }{\partial x^{i}})}V,
 \nabla _{ df_{t}(\f{\partial  }{\partial x^{j}})}V\,\rangle  .
\end{split}\]
Symmetrizing the indexes, the second variation formula becomes
\begin{equation}\label{second_variation}
\begin{split}\f{d^{2} A_{t}}{dt^{2}}&= \int_{M}
(\sum_{i,j}G^{ij} \langle \,\nabla _{df (\f{\partial }{\partial
x^{i}})}V,
 \nabla _{ df (\f{\partial  }{\partial x^{j}})}V\,\rangle
-\f{1}{2}  \sum _{i,j,k,l} G^{ik} \dot{G}_{kl}G^{lj}\dot{G}_{ij})
\,dv_{t}\\&+\int_{M} \sum_{i,j}G^{ij} \langle \,R(V,df(\f{\partial
}{\partial x^{j}}))V, df(\f{\partial }{\partial x^{i}})\,\rangle
\,dv_{t}+\f{1}{4}\int_{M} (\sum_{i,j} G^{ij}\dot{G}_{ij}
)^{2}\,dv_{t}\\
&+\int_{M} \sum_{i,j}G^{ij} \langle \,\nabla _{df(\f{\partial
}{\partial x^{i}})}\nabla _{V}V, df(\f{\partial }{\partial
x^{j}})\,\rangle\,dv_{t}.\\
\end{split}\end{equation}

This formula will be used to prove the main theorems in the next
section.

\section{The stability and uniqueness of minimal maps}

We recall a minimal submanifold is called {\it stable} if the
second derivative of the volume functional with respect to any
compact supported normal variation is non-negative. We prove the
following lemma for minimal graphs.
\begin{lem}
Suppose that the graph of $f:M\rightarrow N$ is a minimal
submanifold $\Sigma$ in $M\times N$.
 Then $\Sigma$ is
stable if and only if it is stable with respect to any compact
supported deformation of maps from $M$ to $N$.
\end{lem}

\begin{proof}
Suppose that
 $a_{i}$ is an
orthonormal basis of the principal  directions of $df$ with
stretches $\lambda_i\geq 0$ and that $df(a_i)=\lambda_{i} b_{i}$.
Assume that the rank of $df(x)$ is $p$.  The orthonormal
set $\{b_{i}\}_{i=1\cdots p}$ can be completed to
form  a local orthonormal basis $\{b_{\alpha}\}_{\alpha=1\cdots
m}$ of the tangent space of $N$.
 In the basis chosen as above, the tangent space of
   $\Sigma$ is spanned by
$t_{i}=\f{1}{\sqrt{1+\lambda_{i}^{2}}}(a_{i}+\lambda_{i}b_{i}),\hskip.2cm
1\leq i \leq n.$ Observe  that  $\lambda_{i}=0$ for $p< i\leq n$.
The normal space of
   $\Sigma$ is spanned by
$n_{i}=\f{1}{\sqrt{1+\lambda_{i}^{2}}}(b_{i}-\lambda_{i}a_{i}),\hskip.2cm
1\leq i \leq p$ and $n_{\alpha}=b_{\alpha}$ for $p<\alpha \leq m$.
Assume that $\bar{V}=\sum_{\alpha=1}^{m}v_{\alpha}n_{\alpha}$ is a
compact supported normal vector field
 along $\Sigma$.  Then the  compact
supported vector field
$V=\sum_{i}\sqrt{1+\lambda_{i}^{2}}\,v_{i}b_{i} +\sum_{\alpha>p}
v_{\alpha}b_{\alpha}$ tangent to $N$ satisfies
$V^{\perp}=\bar{V},$ where $(\cdot)^\perp$ denotes the normal part
of a vector, i.e. the projection onto the normal space of
$\Sigma$.
 The second derivative of volume functional in the direction $V^{\perp}=\bar{V}$ is
  the same as
 in the direction $V$. The Lemma is thus proved.

\end{proof}

The notion of a (strictly) distance-decreasing map in Definition~1
can be generalized to maps between Riemannian manifolds and we can prove
the following theorem.
\begin{thm}
Suppose that $M$ and $N$ are two  Riemannian manifolds, where the
sectional curvature of $N$ is non-positive. Assume that
$f:M\rightarrow N $
  is a  distance-decreasing map and the graph of $f$,
  which is denoted by $\Sigma$, is minimal  in $M\times N$. Then the minimal
  submanifold $\Sigma$
  is stable. It is strictly stable in the following two cases:

  (i) $N$ has negative sectional curvature, and $f$ is not a constant map.

 (ii) $f$ is strictly distance-decreasing, and $M$ is noncompact or with nonempty
 boundary.
 \end{thm}
 \begin{proof}
For a minimal submanifold, we have $\f{dA_t}{dt}|_{t=0}=0$ for any
variation field and in particular
 $$\int_{M}\sum_{i,j}G^{ij}
\langle \,\nabla _{df(\f{\partial  }{\partial x^{i}})}\nabla
_{V}V, df(\f{\partial  }{\partial x^{j}})\,\rangle  \,dv=0.$$ In
the  basis chosen in the proof of Lemma 3.1, we derive from
(\ref{second_variation})
\[\begin{split}\f{d^{2} A_{t}}{dt^{2}}|_{t=0}&\geq\int_{M}(\,\sum_{i}
\f{1}{1+\lambda_{i}^{2}}(| \nabla _{df(a_{i})}V|^{2}-\langle
\,R(V,df(a_{i}))df(a_{i}),V\,\rangle  )
\\&\hskip.5cm-\f{1}{2}\sum_{i,j}
\f{1}{1+\lambda_{i}^{2}}\f{1}{1+\lambda_{j}^{2}}(\langle \,\nabla
_{df(a_{i})}V ,df(a_{j})\,\rangle +\langle  \,\nabla
_{df(a_{j})}V,df(a_{i})\,\rangle  )^{2}\,)\,dv.
\end{split}\]
Since the sectional curvature of $N$ is non-positive, this becomes
\begin{equation}\label{M2nd}
\begin{split}\f{d^{2} A_{t}}{dt^{2}}|_{t=0}&\geq
 \int_{M}(\,\sum_{i}
\f{1}{1+\lambda_{i}^{2}}|\nabla
_{df(a_{i})}V|^{2}\\&\hskip.5cm-\f{1}{2}\sum_{i,j}
\f{1}{1+\lambda_{i}^{2}}\f{1}{1+\lambda_{j}^{2}}(\lambda_{j}
\langle \,\nabla _{df(a_{i})}V,b_{ j}\,\rangle +\lambda_{i}\langle
\,\nabla _{df(a_{j})}V,b_{ i}\,\rangle  )^{2}\,)\,dv
\\&\geq
\int_{M}(\,\sum_{i,j} \f{1}{1+\lambda_{i}^{2}}\langle  \,\nabla
_{df(a_{i})}V,b_{j}\,\rangle  ^{2}\\&\hskip.5cm-\sum_{i,j}
\f{1}{1+\lambda_{i}^{2}}\f{1}{1+\lambda_{j}^{2}}(\lambda_{j}^{2}
\langle  \,\nabla _{df(a_{i})}V,b_{j}\,\rangle  ^{2}
+\lambda_{i}^{2}\langle  \,\nabla _{df(a_{j})}V,b_{i}\,\rangle
^{2})\,)\,dv
\\&=\int_{M}\sum_{i,j}\f{\langle \,\nabla _{df(a_{i})}V,b_{j}\,\rangle  ^{2}}
{1+\lambda_{i}^{2}}\:\f{1-\lambda_{j}^{2}}{1+\lambda_{j}^{2}}\,dv.
\end{split}
\end{equation}

\no When $f$ is a distance-decreasing map, we have $\lambda
_{j}\leq 1$ for $1\leq j\leq n$.  From the estimate in
(\ref{M2nd}), it follows that $\f{d^{2} A_{t}}{dt^{2}}|_{t=0}\geq
0$. This implies that $\Sigma$ is stable by  Lemma 3.1.  Suppose that
$f$ is strictly distance-decreasing, i.e.  $\lambda _{j}< 1$ for
$1\leq j\leq n.$ If $\f{d^{2} A_{t}}{dt^{2}}|_{t=0}= 0$,
 it implies that $\langle \,\nabla _{df(a_{i})}V,b_{j}\,\rangle  =0$  for $1\leq i,j\leq n$ and
  $|\nabla _{df(a_{i})}V|^{2}=\sum _{j}\langle \,\nabla _{df(a_{i})}V,b_{j}\,\rangle  ^{2}$. Hence
    $\nabla _{df(a_{i})}V=0$ for  $1\leq i\leq n$. That is, $V$ is a parallel
    vector field.  In case (ii), $V$ either vanishes outside a compact set or
     on the boundary
    of $M$, so the parallel condition  implies that $V$ is a zero vector.
    This proves that $\Sigma$ is strictly stable in case (ii).
 When the sectional curvature of $N$ is
negative and $f$ is not a constant map, one always has $\f{d^{2} A_{t}}{dt^{2}}|_{t=0}> 0$
unless $V$ is a zero vector.  Therefore,
 $\Sigma$ is strictly stable in case (i).

  \end{proof}
\begin{rem} In case that $M$ is compact without boundary and $f$
is strictly distance-decreasing, one still has the following conclusion:
If $\f{d^{2} A_{t}}{dt^{2}}|_{t=0}= 0$, then
 $V$ is a parallel vector field and $\langle \,R(V,df_{0}(a_{i}))df_{0}(a_{i}),V\,
 \rangle  =0\hskip.3cm$for
$1\leq i\leq n.$
     \end{rem}
  Using the second variation formula, we can also prove the uniqueness of
  minimal maps.
  \begin{thm}
  Suppose that $M$ and $N$ are two  Riemannian manifolds and the sectional
curvature of $N$ is non-positive. Let $\Sigma_{0}$ and
$\Sigma_{1}$ be  minimal
  submanifolds in $M \times N$, which are the graphs
 of distance-decreasing maps $f_{0}:M\rightarrow N$ and
 $f_{1}:M\rightarrow N$, respectively.
 Assume that $f_0$ and $f_1$ are
 homotopic,  and are identical
 on the boundary of $M$ and outside a compact set of $M$.
 Then $\Sigma_{0}=\Sigma_{1}$ in the following two cases:

(i) The sectional curvature of $N$ is negative, and $f_{i}$ are not constant maps,

(ii) The boundary of $M$ is nonempty, or $M$ is noncompact.

\end{thm}
\begin{proof}
 Lift the homotopy map between $f_0$ and $f_1$ to the universal covering of $N$.
Because the sectional curvature of $N$ is non-positive, there exists
a unique geodesic connecting the lifting $\widetilde{f}_{0}(x)$ and
$\widetilde{f}_{1}(x). $ Denote the projection of this unique geodesic
onto $N$  by $\gamma _{x}(t)$ and  define $f_{t}(x)=\gamma _{x}(t)$. Then
$V=\dot{\gamma}_{x}(t)$ satisfies $\nabla _{V}V=0$. Hence the same
bound on $\f{d^{2} A_{t}}{dt^{2}}$ as in (\ref{M2nd}) holds for
$0\leq t\leq 1$.

 The vector field $df_{t}(\f{\partial }{\partial x^{i}})$
  is a Jacobi
field along $\gamma_{x}(t)$, which is denoted by $J_{i,x}(t)$. A
direct calculation gives
\begin{equation} \label{Jacobi}
\begin{split}
\f{d^{2}}{dt^{2}}|J_{i,x} |^{2}&=2\langle  \,\ddot{J}_{i,x}
,J_{i,x}\,
 \rangle+2 |\dot{J}_{i,x} |^{2}
=2\langle  \,R(V,J_{i,x})V,J_{i,x}\,\rangle  +2 |\dot{J}_{i,x}
|^{2}\geq 0.
\end{split}\end{equation}
The last inequality follows from the fact that $N$ has nonpositive
sectional curvature. Because both $f_{0}$ and $f_{1}$ are
distance-decreasing  maps, one has  $|J_{i,x}(0)|^{2}\leq
|\f{\partial }{\partial x^{i}}|^{2}$ and $|J_{i,x}(1)|^{2}\leq
|\f{\partial }{\partial x^{i}}|^{2}$.  The inequality
(\ref{Jacobi}) then implies $|J_{i,x}(t)|^{2}\leq |\f{\partial
}{\partial x^{i}}|^{2}$. Hence $f_{t}$ is also
distance-decreasing and one concludes that  $\f{d^{2}
A_{t}}{dt^{2}}\geq 0$ from (\ref{M2nd}) for $0\leq t\leq 1$.
Because $\f{d  A_{t}}{dt}|_{t=0}=\f{d  A_{t}}{dt}|_{t=1}=0$, the
bound gives $\f{d  A_{t}}{dt}=0$ and $\f{d^{2} A_{t}}{dt^{2}}= 0$
for $0\leq t\leq 1$. To have $\f{d^{2} A_{t}}{dt^{2}}|_{t=0}= 0,$
the following conditions must hold:
\begin{enumerate}
\item $\sum_{i}\f{1}{1+\lambda_{i}^{2}}\langle \,\nabla
_{df_{0}(a_{i})}V,df_{0}(a_{i})\,\rangle  =0.$ \item $  |\nabla
_{df_{0}(a_{i})}V|^{2}=\sum_{j}\langle \,\nabla
_{df_{0}(a_{i})}V,b_{ j}\,\rangle ^{2}\hskip.3cm$for $1\leq i\leq
n.$ \item $  \langle \,\nabla
_{df_{0}(a_{i})}V,df_{0}(a_{j})\,\rangle =\langle \,\nabla
_{df_{0}(a_{j})}V,df_{0}(a_{i})\,\rangle  \hskip.3cm$for $1\leq
i,j\leq n.$ \item If $\lambda_{j}<1$, then $\langle \,\nabla
_{df_{0}(a_{i})}V,b_{j}\,\rangle  =0\hskip.3cm$for $1\leq i\leq
n,$ which implies  $\langle \,\nabla
_{df_{0}(a_{i})}V,df_{0}(a_{j})\,\rangle  =0$. \item $\langle
\,R(V,df_{0}(a_{i}))df_{0}(a_{i}),V\,\rangle  =0\hskip.3cm$for
$1\leq i\leq n.$
\end{enumerate}
When the sectional curvature of $N$ is negative and $f_{0}$ is not a
constant map, condition $5$
implies that $V=0$. Hence $f_{0}=f_{1}$ and
$\Sigma_{0}=\Sigma_{1}$.

Now suppose that the sectional curvature of $N$ is non-positive,
we shall conclude $\nabla _{df_{0}(a_{i})}V=0$ for any $1\leq
i\leq n $. Fix a point $x\in M$ and choose coordinates at $x$ such
that $a_{i}=\f{\partial}{\partial x^{i}}$ for $1\leq i\leq n.$ If
$\lambda_{i}=1$, we have $|df_{0}(\f{\partial}{\partial
x^{i}})|^{2}=1$ and $|J_{i,x}(t)|^{2}$ achieves its maximum at
$t=0$. Therefore,  we have $\f{d }{dt}|J_{i,x}(t)|^{2} =0$ and
$\f{d^{2} }{dt^{2}}|J_{i,x}(t)|^{2}\leq 0$ at $t=0$. The bound on
(\ref{Jacobi}) then implies  $\dot{J}_{i,x}(0)=0$. Hence $\nabla
_{ df_{0}(a_{i})}V =\nabla _{ df_{0}(\f{\partial}{\partial
x^{i}})} V =\nabla _{V}df_{0}(\f{\partial}{\partial x^{i}})=0$. If
$\lambda_{i}<1$,
 condition $3$ and  $4$ give $\langle \,\nabla _{df_{0}(a_{i})}V,df_{0}(a_{j})\,\rangle
 =\langle \,\nabla _{df_{0}(a_{j})}V,df_{0}(a_{i})\,\rangle
 =0$  for $1\leq j\leq n$. Hence $\langle \,\nabla _{df_{0}(a_{i})}V, b_{j} \,\rangle  =0$
 if $\lambda _{j} \neq 0$.  One can still conclude that $\langle \,\nabla
_{df_{0}(a_{i})}V, b_{j} \,\rangle  =0$ from condition $4$ in case $\lambda_j=0$.
 Condition $2$ then implies $\nabla _{df_{0}(a_{i})}V=0$ in the
case $\lambda_{i}<1$.

In conclusion, we always have $\nabla _{df_{0}(a_{i})}V=0$  for any $1\leq i\leq n $ and $V$ is a parallel vector field.
In case (ii), the variation field $V$ either vanishes on the boundary or outside a compact set of $M$. It thus implies  $V=0$ on $  M$. Therefore, $f_{0}=f_{1}$
and $\Sigma_{0}=\Sigma_{1}$ in case (ii).

\end{proof}
\begin{rem}When $M$ is
compact without boundary and $N$ has negative sectional curvature,
then either $f_{0}=f_{1}$ or both $f_{0}$ and $f_{1}$ are constants. If
we only know that $N$ has non-positive sectional curvature, we
 can still conclude that
$V$ is a parallel vector field on $f_{t}(M)$ for $0\leq t\leq 1$.
The graphs of $f_{t},\hskip.2cm 0\leq t\leq 1,$ are then minimal
submanifolds of constant distance. Moreover, the Jacobi fields
$J_{i,x}(t)=df_{t}(\f{\partial}{\partial x^{i}}),\hskip.2cm i=1,\cdots, n$
are parallel along $\gamma_{x}(t)$.
 It implies
that the induced metrics on the graphs of $f_{t}$ are the same.
We also have
$\dot{J}_{i,x}(t)=0$ and $\ddot{J}_{i,x}(t)=0$. The Jacobi
equation thus leads to $R(V,df_{t}(\f{\partial}{\partial
x^{i}}))V=0$ for $1\leq i\leq n $ and $0\leq t\leq 1$.
Hence $\langle R(V,T)V, T \rangle=0$ for any vector $T$ tangent to $f_{t}(M)$
in $N$. The results and further exploration are very similar to the case
of harmonic maps as studied by Schoen and Yau in \cite{sy2}.
\end{rem}

\section{Another criterion for stability}
In this section, we will derive another criterion for the stability of minimal
maps. It is in terms of bounds on the two-Jacobian $|\wedge ^{2} df|(x)$
as defined in the introduction. The theorem generalizes the results for nonparametric
minimal submanifolds of codimension one.
\begin{thm}
Let $M$ and $N$ be Riemannian manifolds and $\Sigma$ be the graph
of a map $f: M\rightarrow N$ with $\text{rank}(df)\leq p$ for some integer $p>1$.  Suppose the sectional curvature of $N$
is non-positive and $\Sigma$ is minimal in $M\times N$. Then
$\Sigma$ is stable if $|\wedge^{2} df|(x)\leq \f{1}{p-1}$ for
any $x\in M$.
 \end{thm}
\begin{proof}
We will keep the term $\f{1}{4}\int_{M} (\sum_{i,j} G^{ij}\dot{G}_{ij}
)^{2}\,dv$ in the second variation formula. In
the  basis chosen in the proof of Lemma 3.1, we derive from
(\ref{second_variation})
\[\begin{split}\f{d^{2} A_{t}}{dt^{2}}|_{t=0}&= \int_{M}(\,\sum_{i}
\f{1}{1+\lambda_{i}^{2}}(| \nabla _{df(a_{i})}V|^{2}-\langle
\,R(V,df(a_{i}))df(a_{i}),V\,\rangle  )
\\&\hskip.5cm-\f{1}{2}\sum_{i,j}
\f{1}{1+\lambda_{i}^{2}}\f{1}{1+\lambda_{j}^{2}}(\langle \,\nabla
_{df(a_{i})}V ,df(a_{j})\,\rangle +\langle  \,\nabla
_{df(a_{j})}V,df(a_{i})\,\rangle  )^{2}\,)\,dv
\\&\hskip.5cm +\int_{M}(\,\sum_{i}
\f{1}{1+\lambda_{i}^{2}}\langle \,\nabla
_{df(a_{i})}V ,df(a_{i})\,\rangle)^{2}\,dv.
\end{split}\]
Since the sectional curvature of $N$ is non-positive, this becomes
\begin{equation}\label{M3rd}
\begin{split}\f{d^{2} A_{t}}{dt^{2}}|_{t=0}&\geq
 \int_{M}(\,\sum_{i}
\f{1}{1+\lambda_{i}^{2}}|\nabla
_{df(a_{i})}V|^{2}+(\,\sum_{i}
\f{\lambda_{i}}{1+\lambda_{i}^{2}}\langle \,\nabla
_{df(a_{i})}V ,b_{i}\,\rangle)^{2}\\&\hskip.5cm-\f{1}{2}\sum_{i,j}
\f{1}{(1+\lambda_{i}^{2})(1+\lambda_{j}^{2})}(\lambda_{j}
\langle \,\nabla _{df(a_{i})}V,b_{ j}\,\rangle +\lambda_{i}\langle
\,\nabla _{df(a_{j})}V,b_{ i}\,\rangle  )^{2}\,)\,dv
\\&\geq
\int_{M}(\,\sum_{i,j} \f{1}{1+\lambda_{i}^{2}}\langle  \,\nabla
_{df(a_{i})}V,b_{j}\,\rangle  ^{2}\\&\hskip.5cm +\sum_{i, j}
\f{\lambda_{i}\lambda_{j}}{(1+\lambda_{i}^{2})(1+\lambda_{j}^{2})}\langle \,\nabla
_{df(a_{i})}V ,b_{i}\,\rangle
\langle \,\nabla
_{df(a_{j})}V ,b_{j}\,\rangle
\\&\hskip.5cm-\f{1}{2}\sum_{i,j}
\f{1}{(1+\lambda_{i}^{2})(1+\lambda_{j}^{2})}(\lambda_{j}
\langle \,\nabla _{df(a_{i})}V,b_{ j}\,\rangle +\lambda_{i}\langle
\,\nabla _{df(a_{j})}V,b_{ i}\,\rangle  )^{2}\,)\,dv
\end{split}
\end{equation}
We break the terms into $i=j$ and $i\neq j$, and obtain
\[\begin{split}
&\sum_{i,j} \f{1}{1+\lambda_{i}^{2}}\langle  \,\nabla
_{df(a_{i})}V,b_{j}\,\rangle  ^{2}\\=&\sum_{i} \f{1}{1+\lambda_{i}^{2}}\langle  \,\nabla
_{df(a_{i})}V,b_{i}\,\rangle  ^{2}+\sum_{i \neq j} \f{1}{1+\lambda_{i}^{2}}\langle  \,\nabla
_{df(a_{i})}V,b_{j}\,\rangle  ^{2},
\end{split}\]
and
\[\begin{split}
&\sum_{i, j}
\f{\lambda_{i}\lambda_{j}}{(1+\lambda_{i}^{2})(1+\lambda_{j}^{2})}\langle \,\nabla
_{df(a_{i})}V ,b_{i}\,\rangle
\langle \,\nabla
_{df(a_{j})}V ,b_{j}\,\rangle\\=&\sum_{i}
\f{\lambda_{i}^{2}}{(1+\lambda_{i}^{2})^{2}}\langle \,\nabla
_{df(a_{i})}V ,b_{i}\,\rangle^{2}+\sum_{i\neq j}
\f{\lambda_{i}\lambda_{j}}{(1+\lambda_{i}^{2})(1+\lambda_{j}^{2})}\langle \,\nabla
_{df(a_{i})}V ,b_{i}\,\rangle
\langle \,\nabla
_{df(a_{j})}V ,b_{j}\,\rangle,
\end{split}\]
and
\[\begin{split}
&\f{1}{2}\sum_{i,j}
\f{1}{(1+\lambda_{i}^{2})(1+\lambda_{j}^{2})}(\lambda_{j}
\langle \,\nabla _{df(a_{i})}V,b_{ j}\,\rangle +\lambda_{i}\langle
\,\nabla _{df(a_{j})}V,b_{ i}\,\rangle  )^{2}
\\=&\sum_{i}\f{2\lambda_{i}^{2}}{(1+\lambda_{i}^{2})^{2}}
\langle  \,\nabla _{df(a_{i})}V,b_{i}\,\rangle  ^{2}
+\sum_{i\neq j}
\f{\lambda_{j}^{2}}{(1+\lambda_{i}^{2})(1+\lambda_{j}^{2})}
\langle  \,\nabla _{df(a_{i})}V,b_{j}\,\rangle  ^{2}\\
&+\sum_{i\neq j}
\f{\lambda_{i}\lambda_{j}}{(1+\lambda_{i}^{2})(1+\lambda_{j}^{2})}\langle  \,\nabla _{df(a_{i})}V,b_{j}\,\rangle
\langle  \,\nabla _{df(a_{j})}V,b_{i}\,\rangle).
\end{split}\]
Plug these expressions into (\ref{M3rd}), and obtain
\begin{equation}\label{M4th}
\begin{split}\f{d^{2} A_{t}}{dt^{2}}|_{t=0}&\geq
\int_{M}(\sum_{i} \f{1}{(1+\lambda_{i}^{2})^{2}}\langle  \,\nabla
_{df(a_{i})}V,b_{i}\,\rangle  ^{2}\\&\hskip.5cm+\sum_{i\neq j}
\f{\lambda_{i}\lambda_{j}}{(1+\lambda_{i}^{2})(1+\lambda_{j}^{2})}\langle \,\nabla
_{df(a_{i})}V ,b_{i}\,\rangle
\langle \,\nabla
_{df(a_{j})}V ,b_{j}\,\rangle\,)\,dv\\&\hskip.5cm+\int_{M}(\sum_{i\neq j}
\f{1}{(1+\lambda_{i}^{2})(1+\lambda_{j}^{2})}\langle  \,\nabla
_{df(a_{i})}V,b_{j}\,\rangle  ^{2}\\&\hskip.5cm-\sum_{i\neq j}
\f{\lambda_{i}\lambda_{j}}{(1+\lambda_{i}^{2})(1+\lambda_{j}^{2})}\langle  \,\nabla _{df(a_{i})}V,b_{j}\,\rangle
\langle  \,\nabla _{df(a_{j})}V,b_{i}\,\rangle
 \,)\,dv.
\end{split}
\end{equation}
The sum of the first two integrands on the right hand side of (\ref{M4th}) is no less than
\[ \begin{split} &\sum_{\lambda_{i}\neq 0}\f{\langle  \,\nabla
_{df(a_{i})}V,b_{i}\,\rangle  ^{2}}{(1+\lambda_{i}^{2})^{2}}+\sum_{i\neq j, \lambda_{i}\neq 0,
\lambda_{j}\neq 0}
\f{\lambda_{i}\lambda_{j}\langle \,\nabla
_{df(a_{i})}V ,b_{i}\,\rangle
\langle \,\nabla
_{df(a_{j})}V ,b_{j}\,\rangle}{(1+\lambda_{i}^{2})(1+\lambda_{j}^{2})}
\\ \geq&\sum_{i\neq j,\lambda_{i}\neq 0,\lambda_{j}\neq 0}\f{\langle  \,\nabla
_{df(a_{i})}V,b_{i}\,\rangle  ^{2}}{(p-1)(1+\lambda_{i}^{2})^{2}}+
\f{\lambda_{i}\lambda_{j}\langle \,\nabla
_{df(a_{i})}V ,b_{i}\,\rangle
\langle \,\nabla
_{df(a_{j})}V ,b_{j}\,\rangle}{(1+\lambda_{i}^{2})(1+\lambda_{j}^{2})}
\\\geq&\sum_{i\neq j,\lambda_{i}\neq 0,\lambda_{j}\neq 0}\f{\langle  \,\nabla
_{df(a_{i})}V,b_{i}\,\rangle  ^{2}}{(p-1)(1+\lambda_{i}^{2})^{2}}
-\f{|\langle \,\nabla
_{df(a_{i})}V ,b_{i}\,\rangle||
\langle \,\nabla
_{df(a_{j})}V ,b_{j}\,\rangle|}{(p-1)(1+\lambda_{i}^{2})(1+\lambda_{j}^{2})}
\\=&\f{1}{p-1}\sum_{i< j,\lambda_{i}\neq 0,\lambda_{j}\neq 0}\f{\langle  \,\nabla
_{df(a_{i})}V,b_{i}\,\rangle  ^{2}}{(1+\lambda_{i}^{2})^{2}}
-2\,\f{|\langle \,\nabla
_{df(a_{i})}V ,b_{i}\,\rangle||
\langle \,\nabla
_{df(a_{j})}V ,b_{j}\,\rangle|}{(1+\lambda_{i}^{2})(1+\lambda_{j}^{2})}
+\f{\langle  \,\nabla
_{df(a_{j})}V,b_{j}\,\rangle  ^{2}}{(1+\lambda_{j}^{2})^{2}}\\=&
\f{1}{p-1}\sum_{i< j,\lambda_{i}\neq 0,\lambda_{j}\neq 0}\,(\,\f{|\langle  \,\nabla
_{df(a_{i})}V,b_{i}\,\rangle| }{1+\lambda_{i}^{2}}-\f{|\langle  \,\nabla
_{df(a_{j})}V,b_{j}\,\rangle| }{1+\lambda_{j}^{2}}
\,)^{2}.
\end{split}\]
While the sum of the last two integrands on the right hand side of (\ref{M4th}), after symmetrizing the indexes, can be written as
$$
\sum_{i\neq j}
\f{\langle  \,\nabla
_{df(a_{i})}V,b_{j}\,\rangle  ^{2}-
2\lambda_{i}\lambda_{j}\langle  \,\nabla _{df(a_{i})}V,b_{j}\,\rangle
\langle  \,\nabla _{df(a_{j})}V,b_{i}\,\rangle+\langle  \,\nabla
_{df(a_{j})}V,b_{i}\,\rangle  ^{2}}{2(1+\lambda_{i}^{2})(1+\lambda_{j}^{2})}.
$$
It is clearly non-negative since $\lambda_{i}\lambda_{j}\leq\f{1}{p-1}\leq 1$
for $i\neq j$. Hence we have $\f{d^{2} A_{t}}{dt^{2}}|_{t=0}\geq 0$ and the
minimal submanifold is stable as claimed.

\end{proof}

\end{document}